\documentclass[11pt]{article}
\usepackage{latexsym,amsmath,color,amsthm,amssymb,epsfig,graphicx,mathrsfs}

\def\qed{\hfill \ifhmode\unskip\nobreak\fi\quad\ifmmode\Box\else$\Box$\fi\\ }
\def\ex{{\rm ex}}

\newtheorem{thrm}{Theorem}[section]
\newtheorem{lemm}[thrm]{Lemma}
\newtheorem{propo}[thrm]{Proposition}
\newtheorem{coro}[thrm]{Corollary}
\newtheorem{defi}[thrm]{Definition}
\newtheorem{conjecture}[thrm]{Conjecture}

\newcommand{\bF}{\mathbb  F}\newcommand{\thm}{\begin{thrm}}
\newcommand{\xthm}{\end{thrm}}
\newcommand{\lem}{\begin{lemm}}
\newcommand{\xlem}{\end{lemm}}
\newcommand{\prf}{\begin{proof}}
\newcommand{\xprf}{\end{proof}}
\newcommand{\prop}{\begin{propo}}
\newcommand{\xprop}{\end{propo}}
\newcommand{\cor}{\begin{coro}}
\newcommand{\xcor}{\end{coro}}
\newcommand{\defn}{\begin{defi}}
\newcommand{\xdefn}{\end{defi}}
\newcommand{\conj}{\begin{conjecture}}
\newcommand{\xconj}{\end{conjecture}}

\renewcommand{\phi}{\varphi}

\newtheorem{question}{Question}

\begin{document}

\title{\vspace{-0.5in}Tur\'an problems and shadows III: expansions of graphs}

\author{
{\large{Alexandr Kostochka}}\thanks{
\footnotesize {University of Illinois at Urbana--Champaign, Urbana, IL 61801
 and Sobolev Institute of Mathematics, Novosibirsk 630090, Russia. E-mail: \texttt {kostochk@math.uiuc.edu}.
 Research of this author
is supported in part by NSF grant  DMS-1266016
and by grants 12-01-00631 and 12-01-00448 
 of the Russian Foundation for Basic Research.
}}
\and
{\large{Dhruv Mubayi}}\thanks{
\footnotesize {Department of Mathematics, Statistics, and Computer Science, University of Illinois at Chicago, Chicago, IL 60607.
E-mail:  \texttt{mubayi@uic.edu.}
Research partially supported by NSF grants DMS-0969092 and DMS-1300138.}}
\and{\large{Jacques Verstra\"ete}}\thanks{Department of Mathematics, University of California at San Diego, 9500
Gilman Drive, La Jolla, California 92093-0112, USA. E-mail: {\tt jverstra@math.ucsd.edu.} Research supported by NSF Grant DMS-1101489. }}
\maketitle

\vspace{-0.3in}

\begin{abstract}
The expansion $G^+$ of a graph $G$ is the $3$-uniform hypergraph obtained from $G$ by enlarging each edge of
$G$ with a new vertex disjoint from $V(G)$ such that distinct edges are enlarged by distinct vertices. Let $\ex_3(n,F)$ denote the maximum number
of edges in a $3$-uniform hypergraph with $n$ vertices not containing any copy of a $3$-uniform hypergraph $F$. The study of $\ex_3(n,G^+)$ includes
some well-researched problems, including the case that $F$ consists of $k$ disjoint edges~\cite{EKR}, $G$ is a triangle~\cite{CK,FF,MV}, $G$ is a path or cycle~\cite{FJS,KMV}, and $G$ is a tree~\cite{ES,Frankl,Furedi,FJ,KMV2}.
In this paper we initiate a broader study of the behavior of $\ex_3(n,G^+)$. Specifically, we show
\[ \ex_3(n,K_{s,t}^+) = \Theta(n^{3 - 3/s})\]
whenever $t > (s - 1)!$ and $s \geq 3$. One of the main open problems is to determine for which graphs $G$ the quantity $\ex_3(n,G^+)$ is quadratic in $n$.
We show that this occurs when $G$ is any bipartite graph with Tur\'{a}n number $o(n^{\varphi})$ where $\varphi = \frac{1 + \sqrt{5}}{2}$,
and in particular, this shows $\ex_3(n,Q^+) = \Theta(n^2)$ where $Q$ is the three-dimensional cube graph.
\end{abstract}

\section{Introduction}

An $r$-uniform hypergraph $F$, or simply {\em $r$-graph}, is a family of $r$-element subsets of a finite set. 
We associate an $r$-graph $F$ with its edge set and call its vertex set $V(F)$. Given an $r$-graph $F$, 
let $\ex_r(n,F)$ denote the maximum number of edges in an $r$-graph on $n$ vertices that does not contain $F$.  
The {\em expansion} of a graph $G$ is the $3$-graph $G^+$ with edge set $\{e \cup \{v_e\} : e \in G\}$ where
$v_e$ are distinct vertices  not in $V(G)$. By definition, the expansion of $G$ has exactly $|G|$ edges.
Note that F\" uredi and Jiang~\cite{Furedi,FJ} used a  notion of expansion to $r$-graphs for general $r$, but 
this paper considers only $3$-graphs.

\medskip

Expansions include many important hypergraphs whose extremal functions have been investigated. For instance, 
the celebrated Erd\H{o}s-Ko-Rado Theorem~\cite{EKR}
for $3$-graphs
 is the case of expansion of a matching. A well-known result is that 
  ex$_3(n, K_3^+)={n-1 \choose 2}$~\cite{CK,FF,MV}. If a graph is not $3$-colorable then its expansion has positive Tur\'{a}n density and this case is fairly well understood~\cite{M,P2}, so we focus on the case of expansions of 3-colorable graphs.
It is easy to see that $\ex_3(n,G^+) = \Omega(n^2)$ unless $G$ is a star (the case that
$G$ is a star is interesting in itself, and for $G = P_2$ determining $\ex_3(n,G^+)$ constituted a conjecture of Erd\H{o}s and S\'{o}s~\cite{ES} which was solved
by Frankl~\cite{Frankl}). The authors~\cite{KMV} had previously determined $\ex_3(n,G^+)$ exactly (for large $n$) when $G$ is a path or cycle of fixed length $k\ge 3$, thereby answering questions of F\"uredi-Jiang-Siever~\cite{FJS} and F\"uredi-Jiang~\cite{FJ}.
The case when $G$ is a forest is solved asymptotically in~\cite{KMV2}, thus settling a conjecture of F\"uredi~\cite{Furedi}. The following straightforward result provides general bounds for $\ex_3(n, G^+)$ in terms of the number of edges of $G$.

\prop \label{planar}
If $G$ is any graph with $v$ vertices and $f\ge 4$ edges, then for some $a > 0$,
\[ an^{3 - \frac{3v - 9}{f - 3}} \leq \ex_3(n,G^+) \leq (n-1)\ex_2(n,G) + (f + v - 1){n \choose 2}.\]
\xprop

The proof of Proposition \ref{planar} is given in Section \ref{planarproof}. Some key remarks are that $\ex_3(n,G^+)$ is not
quadratic in $n$ if $f > 3v - 6$, and if $G$ is not bipartite then the upper bound in Proposition \ref{planar} is cubic in $n$.
This suggests the question of identifying the graphs $G$ for which ex$_3(n,G^+)=O(n^2)$, and in particular
evaluation of $\ex_3(n,G^+)$ for planar $G$.

\subsection{Expansions of planar graphs}

We give a straightforward proof of the following proposition, which is a
special case of a more general result of F\"{u}redi~\cite{Furedi} for a larger class of triple systems.

\prop \label{tw2}
Let $G$ be a graph with treewidth at most two. Then $\ex_3(n,G^+) = O(n^2)$.
\xprop

On the other hand, there are 3-colorable planar graphs $G$ for which $\ex_3(n,G^+)$ is not quadratic in $n$.
 To state this result,
we need a definition. A proper $k$-coloring $\chi : V(G) \rightarrow \{1,\ldots,k\}$ is  {\em acyclic }
if every pair of color classes induces a forest in $G$. We pose the following question:

\begin{question} \label{problem1}
Does every planar graph $G$ with an acyclic 3-coloring have {\rm ex}$_3(n, G^+) = O(n^2)$?
\end{question}

Let $g(n,k)$ denote the maximum number of edges in an $n$-vertex graph of girth larger than $k$.

\prop \label{novalid}
Let $G$ be a planar graph such that in every proper 3-coloring of $G$, every pair of color classes induces a
subgraph containing a cycle of length at most $k$. Then $\ex_3(n,G^+) = \Omega(ng(n,k)) = \Omega( n^{2 + \Theta(\frac{1}{k})})$.
\xprop

The last statement follows from the known fact that $g(n,k) \geq n^{1 + \Theta(\frac{1}{k})}$.
The octahedron graph $O$ is an example of a planar graph where in every proper 3-coloring, each pair of color classes
induces a cycle of length four, and so $\ex_3(n,O^+) = \Omega(n^{5/2})$. Even wheels do not have acyclic 3-colorings,
and we do not know whether their expansions have quadratic Tur\'{a}n numbers.

\begin{question} \label{evenwheel}
Does every even wheel $G$ have $\ex_3(n,G^+) = O(n^2)$?
\end{question}

\subsection{Expansions of bipartite graphs}

The behavior of $\ex_3(n,G^+)$ when $G$ is a dense bipartite graph is somewhat related to the behavior of $\ex_2(n,G)$ according to Proposition \ref{planar}.
In particular, Proposition \ref{planar} shows that for $t \geq s \geq 2$ and some constants $a,c > 0$, 
\[ an^{3 - \frac{3s + 3t - 9}{st - 3}} \leq \ex_3(n,K_{s,t}^+) \leq c n^{3 - \frac{1}{s}}. \]
We show that both the upper and lower bound can be improved to determine the order of magnitude of $\ex_3(n,K_{s,t}^+)$ when
good constructions of $K_{s,t}$-free graphs are available (see Alon, R\'{o}nyai and Szabo~\cite{ARS}):

\thm \label{kst}
Fix $3 \le s \le t$. Then $\ex_3(n,K_{s,t}^+) = O(n^{3 - \frac{3}{s}})$ and, if  $t > (s - 1)! \geq 2$, then $\ex_3(n,K_{s,t}^+) = \Theta(n^{3- \frac{3}{s}})$.
\xthm

The case of $K_{3,t}$ is interesting since $\ex_3(n,K_{3,t}^+) = O(n^2)$, and perhaps it is possible to determine a constant $c$ such that $\ex_3(n,K_{3,3}^+) \sim cn^2$,
since the asymptotic behavior of $\ex_2(n,K_{3,3})$ is known, due to a construction of Brown~\cite{Brown} and the upper bounds of F\"{u}redi~\cite{Furedi}.
 In general, the following bounds hold for expansions of $K_{3,t}$:

\thm \label{k3t}
For fixed $r \ge 1$ and $t = 2r^2 + 1$, we have  $(1-o(1))\frac{t - 1}{12} n^2 \leq \ex_3(n,K_{3,t}^+) = O(n^2)$.
\xthm

The upper bound in this theorem is a special case of a general upper bound for all graphs $G$ with $\sigma(G^+) = 3$ (see Theorem \ref{sigma3}).
Finally, we prove a general result that applies to expansions of a large class of bipartite graphs.

\thm \label{cube}
Let $G$ be a graph with $\ex_2(n, G) = o(n^{\phi})$, where $\phi=(1+\sqrt 5)/2$ is the golden ratio.  Then $\ex_3(n, G^+)= O(n^2)$.
\xthm

Let $\mathbb Q$ be the graph of the 3-dimensional cube (with $8$ vertices and $12$ edges). Erd\H{o}s and Simonovits~\cite{ES} proved $\ex_2(n,\mathbb Q) = O(n^{1.6}) = o(n^{\phi})$, so  a corollary to Theorem \ref{cube} is that
$$\ex_3(n, \mathbb Q^+) = \Theta(n^2).$$
 Determining the growth rate of $\ex_2(n, \mathbb Q)$ is a longstanding open problem. Since it is known that for any graph $G$ the 1-subdivision
of $G$ has Tur\'{a}n Number $O(n^{3/2})$ -- see Alon, Krivelevich and Sudakov~\cite{AKS} -- Theorem \ref{cube} also shows that for such graphs $G$, $\ex_3(n,G^+) = \Theta(n^2)$.
Erd\H{o}s conjectured that  $\ex_2(n,G) = O(n^{3/2})$ for each 2-degenerate bipartite graph $G$. If this conjecture is true, then by Theorem \ref{cube}, 
$\ex_3(n,G^+) = O(n^2)$ for any 2-degenerate bipartite graph $G$.

\subsection{Crosscuts}

A set of vertices in a hypergraph containing exactly one vertex from every edge of a hypergraph is called a {\em crosscut} of the hypergraph,
following Frankl and F\"{u}redi~\cite{FF}.
For a 3-uniform  hypergraph $F$, let $\sigma(F)$ be the minimum size of a crosscut of $F$ if it exists, i.e.,
\[ \sigma(F) := \min\{|X| : \forall e \in F, |e \cap X| = 1\}\]
if such an $X$ exists.
Since the triple system consisting of all edges containing exactly one vertex from a set of size $\sigma(F) - 1$
does not contain $F$, we have
\begin{equation} \label{sigma} \ex_3(n,F) \geq (\sigma(F) - 1 +o(1)){n \choose 2}.
\end{equation}
An intriguing open question is: For which $F$ an asymptotic equality is attained in~\eqref{sigma}? 
Recall that a graph has tree-width at most two if and only if it has no subdivision of $K_4$. Informally, these are subgraphs of a planar graph
obtained by starting with a triangle, and then picking some edge $uv$ of the current graph, adding a new vertex $w$, and then adding the edges $uw$ and $vw$.

\begin{question} \label{problem2}
Is it true that
\begin{equation}\label{q3}
 {\rm ex}_3(n, G^+) \sim (\sigma(G^+) - 1){n \choose 2}
 \end{equation}
 for every graph $G$ with tree-width two?
\end{question}

If $G$ is a forest or a cycle, then \eqref{q3} holds~\cite{KMV,KMV2} (corresponding results for $r>3$ were given by F\"uredi~\cite{Furedi}). 
 If $G$ is any graph with $\sigma(G^+) = 2$,
then again~ \eqref{q3} holds~\cite{KMV2}.  Proposition \ref{planar} and Theorem \ref{kst}
give examples of graphs $G$ with $\sigma(G^+) = 4$ and $\ex_3(n,G^+)$ superquadratic in $n$.
This leaves the case $\sigma(G^+) = 3$,
and in this case, Theorem \ref{k3t} shows that $\ex_3(n,K_{3,t}^+)/n^2 \rightarrow \infty$ as $t \rightarrow \infty$, even though
$\sigma(K_{3,t}^+) = 3$ for all $t \geq 3$. A quadratic upper bound for $\ex_3(n,K_{3,t}^+)$ in Theorem \ref{k3t} is a special case of the following theorem:

\thm \label{sigma3}
For every $G$ with $\sigma(G^+)=3$, $\ex_3(n, G^+)= O(n^2)$.
\xthm


\section{Preliminaries}

{\bf Notation and terminology.} A $3$-graph is called a {\em triple system}. The edges will be written
as unordered lists, for instance, $xyz$ represents $\{x,y,z\}$.
For a set $X$ of vertices of a hypergraph $H$, let $H - X = \{e \in H : e \cap X = \emptyset\}$.  If $X = \{x\}$, then we write
$H - x$ instead of $H - X$.
The {\em codegree} of a pair $\{x,y\}$ of vertices in $H$ is $d_H(x,y) = |\{e \in H : S \subset e\}|$
and for a set $S$ of vertices, $N_H(S)=\{x \in V(H) : S \cup \{x\} \in H\}$ so that $|N_H(S)|=d_H(S)$ when $|S| = 2$.
The {\em shadow} of $H$ is the graph $\partial H = \{xy : \exists \,e \in H, \{x,y\} \subset e\}$.
The edges of $\partial H$ will be called the {\em sub-edges} of $H$.





An $r$-graph  $H$ is {\em $d$-full} if every sub-edge of $H$ has codegree at least $d$.

Thus $H$ is $d$-full is equivalent to the fact that the minimum non-zero codegree in $H$ is at least $d$.
The following lemma  from~\cite{KMV2}
extends the well-known fact that each graph $G$ has a subgraph of minimum degree at least $d$ with at least $|G| - (d-1)|V(G)|$ edges.

\lem \label{fullsub}
For $r \ge 2, d \geq 1$, every $n$-vertex $r$-graph $H$ has a $(d + 1)$-full subgraph $F$ with
\[ |F| \geq |H| - d\,|\partial H|.\]
\xlem

\prf
A {\em $d$-sparse sequence} is a maximal sequence $e_1,e_2,\dots,e_m \in \partial H$ such that $d_H(e_1) \leq d$, and for all $i > 1$, $e_i$
 is contained in at most $d$ edges of $H$ which contain none of $e_1,e_2,\dots,e_{i - 1}$.
The $r$-graph $F$ obtained by deleting all edges of $H$ containing at least one of the $e_i$ is $(d + 1)$-full. Since
a $d$-sparse sequence has length at most $|\partial H|$, we have $|F| \geq |H| - d|\partial H|$.
\xprf

\section{Proofs of Propositions}\label{planarproof}

{\bf Proof of Proposition \ref{planar}.} The proof of the lower bound in Proposition \ref{planar} is via a random triple system. The idea is to take a random graph not containing a particular graph $G$,
and then observe that the triple system of triangles in the random graph does not contain $G^+$.
Consider the random graph on $n$ vertices, whose edges are placed independently with probability $p$, to be chosen later. If $X$ is the number of triangles
and $Y$ is the number of copies of $G$ in the random graph, then
\[ \mathbb E(X) = p^3 {n \choose 3} \quad \quad \mathbb E(Y) \leq p^{f} n^{v}.\]
Therefore choosing $p = 0.1 n^{-(v - 3)/(f - 3)}$, since $f \geq 4$, we find
\[ \mathbb E(X - Y) \geq 0.0001 n^{3 - 3(v-3)/(f - 3)}.\]
Now let $H$ be the triple system of vertex-sets of triangles in the graph obtained by removing one edge from each copy of
$G$ in the random graph. Then $\mathbb E(|H|) \geq \mathbb E(X - Y)$, and $G^+ \not \subset H$.
Select an $H$ so that $|H| \geq 0.0001 n^{3 - 3(v - 3)/(f - 3)}$. This proves the lower bound in Proposition \ref{planar} with $a = 0.0001$.

\medskip

Now suppose $G$ is a bipartite graph with $f$ edges $e_1,e_2,\dots,e_f$ and $v$ vertices. If a triple system $H$ on $n$ vertices has more than $(n - 1)\ex_2(n,G) + (f + v - 1){n \choose 2}$ triples, then by deleting at most $(f + v - 1){n  \choose 2}$ triples we arrive at a triple system $H' \subset H$ which is $(f+v)$-full, by Lemma \ref{fullsub} and $|H'| > (n - 1)\ex_2(n,G)$.
There exists $x \in V(H')$ such that more than $\ex_2(n,G)$ triples of $H'$ contain $x$. So the graph of all pairs $\{x,y\}$ such that $\{w,x,y\} \in H'$
contains $G$. Since every pair $\{w,y\}$ has codegree at least $f + v$, we find vertices $z_1,z_2,\dots,z_f \not \in V(G)$ such
$e_i \cup \{z_i\} \in H'$ for all $i = 1,2,\dots,f$, and this forms a copy of $G^+$ in $H'$.  \qed

\bigskip

{\bf Proof of Proposition \ref{tw2}.}
Let $G$ be a graph of tree-width two. Then $G \subset F$, where $F$ is a graph obtained from a triangle by repeatedly
adding a new vertex and joining it to two adjacent vertices of the current graph. It is enough to show $\ex_3(n,F^+) = O(n^2)$.
Suppose $F$ has $v$ vertices and $f$ edges. By definition, $F$ has a vertex $x$ of degree two such that the neighbors
$x'$ and $x''$ of $x$ are adjacent. Then $F':=F-x$
has $v - 1$ vertices and $f - 2$ edges. Let $H$ be an $n$-vertex triple system with more than $(v + f - 1){n \choose 2}$ edges.
By Lemma \ref{fullsub}, $H$ has a $(v + f)$-full subgraph $H'$. We claim $H'$ contains $F^+$. Inductively, $H'$ contains a copy
$H''$ of the expansion
of $F'$. By the definition of $H'$, $\{x',x''\}$ has codegree at least $v + f$ in $H'$. Therefore we may select a new vertex
$z$ that is not in $H''$ such that $\{z,x',x''\}$ is an edge of $H'$, and now $F$ is embedded in $H'$ by mapping $x$ to $z$. \qed

\bigskip

{\bf Proof of Proposition \ref{novalid}.} Let $G$ be a 3-colorable planar graph with the given conditions.
To show $\ex_3(n,G^+) = \Omega(n g(n,k))$, form a triple system $H$ on $n$ vertices as follows. Let $F$
be a bipartite  $\lfloor \frac{n}{2} \rfloor$-vertex graph of girth $k+1$ with at least $\frac{1}{2}g(\lfloor \frac{n}{2} \rfloor,k)$
edges. Let 
 $U$ and $V$ be the partite sets of $F$.  Let $X$ be a set of $\lceil \frac{n}{2} \rceil$
vertices disjoint from $U\cup V$. Then set $V(H) = U \cup V \cup X$ and let the edges of $H$ consist of all triples $e \cup \{x\}$ such that
$e \in F$ and $x \in X$. Then
\[ |H| \geq |X| \cdot g(\lfloor \frac{n}{2} \rfloor,k) = \Omega(ng(n,k)).\]
Now $\partial H$ has a natural 3-coloring given by $U,V,X$. If $G^+ \subset H$, then $G \subset \partial H$ and therefore
$G$ is properly colored, with color classes $V(G) \cap U$, $V(G) \cap V$ and $V(G) \cap X$. 
By the assumptions on $G$, $V(G) \cap (U\cup V)$ 
induces a subgraph of $G$ which  contains a cycle of length at most $k$. However, that cycle is then
a subgraph of $F$, by the definition of $H$, which is a contradiction. Therefore $G^+ \not \subset H$. \qed

\section{Proof of Theorem \ref{cube}}

{\bf Proof of Theorem \ref{cube}.}
Suppose $\ex_2(n,G) = o(n^{\phi})$ and $|G| = k$, and $H$ is an $G^+$-free $3$-graph
with $|H| \geq (k + 1) {n \choose 2}$. By Lemma \ref{fullsub}, $H$ has a $k$-full-subgraph
$H_1$ with at least $ n^2/3$ edges. If
$G \subset \partial H_1$, then we can expand $G$ to $G^+ \subset H_1$ using that $H_1$ is $k$-full.
Therefore $|\partial H_1| \leq  \ex_2(n,G)=o(n^{\phi})$. By Lemma \ref{fullsub}, and since $|H_1| \geq \delta n^2$,
$H_1$ has a non-empty $n^{2 - \phi}$-full subgraph $H_2$ if $n$ is large enough. Let $H_3$ be obtained
by removing all isolated vertices of $H_2$ and let $m = |V(H_3)|$. Since $H_3$ is
$n^{2 - \phi}$-full, $m> n^{2 - \phi}$.
Since $H_1$ is $G^+$-free,
$H_3 \subset H_1$ is also $G^+$-free, and therefore if $F = \partial H_3$, $|V(F)| = |V(H_3)| = m$ and
$|F| \leq \ex_2(m,G) = o(m^{\phi})$.  So some vertex $v$ of the graph $F = \partial H_3$
has degree $o(m^{\phi - 1})$. Now the number of edges of $F$ between vertices of $N_F(v)$ is at least the number
of edges of $H_3$ containing $v$. Since $H_3$ is $n^{2 - \phi}$-full, there are at least $\frac{1}{2}n^{2 - \phi}|N_F(v)|$
such edges. On the other hand, since the subgraph of $F$ induced by $N_F(v)$ does not contain $G$, the number of such
edges is
 $ o(|N_F(v)|^{\phi})$. It follows that $n^{2 - \phi}=o(|N_F(v)|^{\phi-1})$.
 Since $|N_F(v)| = o(m^{\phi - 1})= o(n^{\phi - 1})$, we get $2-\phi<(\phi-1)^2$, contradicting the fact that
   $\phi$ is the golden ratio. \qed

\section{Proof of Theorem \ref{kst}}

{\bf Proof of Theorem \ref{kst}.} For the upper bound, we repeat the proof of Theorem \ref{cube} when $F = K_{s,t}$,
using the bounds $\ex_2(n,K_{s,t}) = O(n^{2 - 1/s})$ provided by the K\"{o}vari-S\'{o}s-Tur\'{a}n Theorem~\cite{KST},
except
at the stage of the proof where we use the bound on $\ex_2(|N_G(v)|,F)$, we may now use
\[ \ex_2(|N_G(v)|,K_{s - 1,t}) = O(|N_G(v)|^{2 - 1/(s - 1)})\]
for if the subgraph of $G$ of edges between $N_G(v)$ contains $K_{s-1,t}$, then by adding $v$ we see $G$ contains $K_{s,t}$.
A calculation gives $|H| = O(n^{3 - 3/s})$.

For the lower bound  we must show that $\ex_3(n,K_{s,t}^+)= \Omega(n^{3 - 3/s})$ if $t > (s - 1)!$.
We will use the {\em projective norm graphs}
defined by Alon, R\'{o}nyai and Szabo~\cite{ARS}.
Given a finite field $\bF_q$ and an integer $s \ge 2$, the norm  is the map $N: \bF_{q^{s-1}}^* \rightarrow \bF_q^*$ given by $N(X)=X^{1+q+\cdots +q^{s-2}}$.  The norm is a (multiplicative) group homomorphism and is the identity map on elements of $\bF_q^*$.  This implies that for each $x \in \bF_q^*$, the number of preimages of $x$ is exactly
\begin{equation} \label{manytoone} \frac{q^{s-1}-1}{q-1}=1+q+\cdots +q^{s-2}.\end{equation}

\defn  Let $q$ be a prime power and $s \ge 2$ be an integer. The projective norm graph $PG(q,s)$ has vertex set $V=\mathbb F_{q^{s-1}} \times \bF_q^*$ and edge set
$$\{(A, b)(B, b): N(A+B)=ab\}.$$
\xdefn

\lem \label{nnl} Fix an integer $s \ge 3$ and a prime power $q$. Let $x \in \bF_q^*$, and $A, B \in \bF_{q^{s-1}}$ with $A \ne B$.  Then the number of $C \in \bF_{q^{s-1}}$ with
\begin{equation} \label{manyC}N\left(\frac{A+C}{ B+C}\right)=x\end{equation}
is at least $q^{s-2}$.
\xlem
\prf
By (\ref{manytoone}) there exist distinct $X_1, \ldots, X_{q^{s-2}+1} \in \bF_{q^{s-1}}^*$ such that $N(X_i)=x$ for each $i$.  As long as $X_i \ne 1$, define
$$C_i=\frac{BX_i-A}{1-X_i}.$$
Then $(A+C_i)/(B+C_i)=X_i$, and $C_i \ne C_j$ for $i \ne j$ since $A \ne B$.
\xprf

\lem \label{tri} Fix an integer $s \ge 3$ and a prime power $q$. The number of triangles in $PG(q,s)$ is at least
$(1-o(1))q^{3s-3}/6$ as $q \rightarrow \infty$.
\xlem
\prf  Pick a vertex $(A,a)$ and then one of its neighbors $(B, b)$.  The number of ways to do this is at least $q^{s-1}(q-1)(q^{s-1}-1)$. Let $x=a/b$ and apply Lemma~\ref{nnl} to obtain at least $q^{s-2}-2$ distinct $C  \not\in \{-A, -B\}$ satisfying (\ref{manyC}).  For each such $C$, define
$$c=\frac{N(A+C)}{a} = \frac{N(B+C)}{b}.$$
Then $(C,c)$ is adjacent to both $(A, a)$ and $(B, b)$.   Each triangle is counted six times in this way and the result follows.
\xprf

For appropriate $n$
the $n$-vertex norm graphs $PG(q,s)$  (for fixed $s$ and large $q$)
have $\Theta(n^{2 - 1/s})$ edges and no $K_{s,t}$. By Lemma~\ref{tri}
the number of triangles in $PG(q,s)$ is $\Theta(n^{3 - 3/s})$. The hypergraph $H$ whose edges are the vertex sets of triangles in $PG(q,s)$ is a $3$-graph
with $\Theta(n^{3 - 3/s})$ edges and no $K_{s,t}^+$. This completes the proof of Theorem \ref{kst}. \qed

\section{Proof of Theorems \ref{k3t} and \ref{sigma3}}

We need the following result.

\thm \label{hypshadow}
Let $F$ be a 3-uniform hypergraph with $v$ vertices and {\rm ex}$_3(n, F)< c{n\choose 2}$.  Then ex$_3(n, (\partial F)^+) < (c+v+|F|){n \choose 2}$.
\xthm

\prf  Suppose we have an $n$ vertex  3-uniform hypergraph $H$ with $|H|> (c+v+|H|){n \choose 2}$. Apply Lemma \ref{fullsub} to obtain a  subhypergraph $H' \subset H$ that is $(v+|F|)$-full with  $|H'|> c{n\choose 2}$. By definition, we may find a copy of $F \subset H'$ and hence a copy of $\partial F \subset \partial H'$.  Because $H'$ is $(v+|F|)$-full, we may expand this copy of $\partial F$ to a copy of $(\partial F)^+ \subset H' \subset H$ as desired.  \xprf

Define $H_t$ to be the 3-uniform hypergraph with vertex set
$\{a, b, x_1, y_1, \ldots, x_t, y_t\}$ and $2t$ edges $x_iy_ia$ and $x_iy_ib$ for all $i \in [t]$.   It is convenient (though not necessary) for us to use the following theorem of the
authors~\cite{MV0}.

\thm \label{jd} {\bf(\cite{MV0})} For each $t \ge 2$, we have
{\rm ex}$_3(n, H_t)< t^4{n \choose 2}$.
\xthm

{\bf Proof of Theorems \ref{k3t} and \ref{sigma3}.}
First we  prove the upper bound in Theorem~\ref{sigma3}. Suppose $\sigma(G^+) \le 3$. This means that $G$ has an independent set $I$ and set $R$ of edges such
that $I$ intersects each edge in $G-R$, and $|I|+|R|\leq 3$.
It follows that $G$ is a subgraph of one of the following graphs (Cases (i) and (ii) correspond to $|I|=1$, 
Case (iii) corresponds to  $|I|=2$, and Case (iv) corresponds to  $|I|=3$):

(i) $K_4-e$ together with a star centered at one of the degree 3  (in $K_4-e$) vertices,

(ii) two triangles sharing a vertex $x$ and a star centered at $x$,

(iii) the graph obtained from $K_{2,t}$ by adding an edge joining two vertices in the part of size $t$,

(iv) $K_{3,t}$.

Now suppose we have a 3-uniform $H$ with $|H|> cn^2$ for some $c>|G|+|V(G)|$. 
Applying Lemma \ref{fullsub}, we find a $c$-full  $H' \subset H$ with $|H'|> c{n \choose 2}$. 
As in the proof of Theorem~\ref{hypshadow}, it is enough to find $G$ in $\partial H'$. Since
$|H'|> c{n \choose 2}$, the codegree of some pair $\{x,y\}$ is at least $c+1$. Then the shadow of the set of
triples in $H'$ containing $\{x,y\}$ contains the graph of the form i). Similarly, $H'$ contains two edges sharing
exactly one vertex, say $x$, and the shadow of the set of
triples in $H'$ containing $x$ contains the graph of the form ii). 
 If $G$ is of the form in iii), we  apply Theorems \ref{hypshadow} and 
\ref{jd} and observe that $\partial H_t \supset G$.  Finally, if $G\subseteq K_{3,t}$ then we apply Theorem \ref{kst}.

For the lower bound in Theorem~\ref{k3t}, we use a slight modification of the construction in Theorem \ref{kst}. Set $s=3$ and let $r|q-1$. Let $Q_r$ denote a subgroup of
$\bF_q^*$ of order $r$.  Define the graph $H=H_r(q)$ with
$V(H)=\bF_{q^{2}} \times \bF_q^*/Q_r$ and two vertices $(A, aQ_r)$ and $(B, bQ_r)$ are adjacent in $H$ if
$N(A+B)\in abQ_r$.  Then $H$ has $n=(q^3-q^{2})/r$ vertices and each vertex has degree $q^{2}-1$. It also follows from~\cite{ARS} that
$H$ has no $K_{3,t}$ where $t=2r^{2}+1$. Now we construct a 3-uniform hypergraph $H'$ with $V(H')=V(H)$ and whose edges are the triangles of $H$.  We must count the number of triangles in $H$ to determine $|H'|$.  For every choice of $(A, a), (B, b)$ in $\bF_{q^2} \times \bF_q^*$, the number of $(C,c) \in  \bF_{q^2} \times \bF_q^*$ with $C \ne A,B$, $N(A+C)=ac$ and $N(B+C)=bc$ is at least $q-2$ by (the proof of) Lemma \ref{tri}. Consequently, the number of $(C,c)$ such that $N(A+C) \in acQ_r$ and $N(B+C)\in bcQ_r$ is at least $r^2(q-2)$. Since $(C,c)$ satisfies these equations iff $(C, cq)$ satisfies these equations for all $q \in Q_r$ (i.e. the solutions come in equivalence classes of size $r$), the number of common neighbors of
$(A, aQ_r)$ and $(B, bQ_r)$ is at least $r(q-2)$.  The number of edges in $H$ is at least $(1-o(1))q^5/2r$, so the number of triangles in $H$ is at least $(1-o(1))q^6/6=(1-o(1))(r^2/6)n^2$.\qed

\section{Concluding remarks}

$\bullet$ In this paper we studied $\ex_3(n,G^+)$ where $G$ is a 3-colorable graph. If $G$ has treewidth two, then we believe $\ex_3(n,G^+) \sim (\sigma(G^+) - 1){n \choose 2}$ (Question~3), and if $G$ has an acyclic 3-coloring, then we believe $\ex_3(n,G^+) = O(n^2)$ (Question~1). We are also not able to prove
or disprove $\ex_3(n,G^+) = O(n^2)$ when $G$ is an even wheel (Question~2). This is equivalent to showing that if $F$ is an $n$-vertex graph with a
superquadratic number of triangles, then $F$ contains every even wheel with a bounded number of vertices.

\medskip

$\bullet$ A number of examples of $3$-colorable $G$ with superquadratic $\ex_3(n,G^+)$ were given. In particular we determined the order of magnitude of $\ex_3(n,K_{s,t}^+)$
when near-extremal constructions of $K_{s,t}$-free bipartite graphs are known. One may ask for the asymptotic behavior of $\ex_3(n,K_{3,t}^+)$ for each $t \geq 3$, since in
that case we have shown $\ex_3(n,K_{3,t}^+) = \Theta(n^2)$.
Finally, we gave a general upper bound on $\ex_3(n,G^+)$ when $G$ is a bipartite graph, and showed that if $G$ has Tur\'{a}n number much smaller than $n^{\varphi}$ where
$\varphi$ is the golden ratio, then $\ex_3(n,G^+) = O(n^2)$. Determining exactly when $\ex_3(n,G^+)$ is quadratic in $n$ remains an open problem for further research.


\begin{thebibliography}{99}

\bibitem{AKS} Alon, N.; Krivelevich, M.; Sudakov, B. Tur\'{a}n numbers of bipartite graphs and related Ramsey-type questions.
Combinatorics, Probability and Computing 12 (2003), 477--494.

\bibitem{ARS} N. Alon, L. R\'{o}nyai and T. Szab\'{o}. Norm-graphs: variations and applications, J. Combinatorial Theory, Ser. B 76 (1999), 280--290.


\bibitem{Brown} W. Brown, On graphs that do not contain a Thomsen graph, Canadian Mathematical
Bulletin 9 (1966), 281--285.


\bibitem{CK}  Cs\'ak\'any, R.; Kahn, J. A homological approach to two problems on finite sets. J. Algebraic Combin. 9 (1999), no. 2, 141--149.




\bibitem{EKR} Erd\H{o}s, P.; Ko, C.; Rado, R. Intersection theorems for systems of finite sets, The Quarterly Journal of Mathematics. Oxford. Second Series (1961) 12, 313--320.


\bibitem{ES} Erd\H{o}s, P. Extremal problems in graph theory in:  Theory of Graphs and its Applications,
M. Fiedler (Ed.), Academic Press, New York (1965), pp. 29--36.


\bibitem{Frankl} Frankl, P. On families of finite sets no two of which intersect in a singleton. Bull. Austral. Math. Soc.
17 (1977), no. 1, 125--134.

\bibitem{FF} Frankl, P.; F\" uredi, Z. Exact solution of some Tur\'{a}n-type problems. J. Combin. Theory Ser. A 45 (1987), no. 2, 226--262.

\bibitem{Furedi}  F\"{u}redi, Z.
 Linear trees in uniform hypergraphs. European J. Combin. Theory Ser. A 35 (2014), 264--272.

\bibitem{FJ} F\"{u}redi, Z.; Jiang, T. Hypergraph
Tur\'an numbers of linear cycles. Preprint (2013). arXiv:1302.2387

\bibitem{FJS} F\"{u}redi, Z.; Jiang, T.; Seiver, R. Exact solution of the hypergraph Tur\'{a}n problem for $k$-uniform linear paths, To appear in Combinatorica (2013).


\bibitem{KMV}  Kostochka, A.; Mubayi, D.; Verstra\"ete, J.; Tur\'an problems and shadows I: paths and cycles, submitted.

\bibitem{KMV2}  Kostochka, A.; Mubayi, D.; Verstra\"ete, J.; Tur\'an problems and shadows  II: trees, submitted.


\bibitem{KST} K\"{o}vari, T.; S\'{o}s, V. T.; Tur\'{a}n, P.
On a problem of K. Zarankiewicz. Colloquium Math. 3, (1954). 50--57.

\bibitem{M}  Mubayi, D. A hypergraph extension of Tur\'{a}n's theorem.
J. of Combinatorial Theory, Ser. B, 96 (2006), no. 1, 122--134.

\bibitem{MV0} Mubayi, D.; Verstra\"{e}te, J. A hypergraph extension of the Bipartite Tur\' an problem, Journal of Combinatorial Theory, Series A 106 (2004) no. 2, 237--253.

\bibitem{MV} Mubayi, D.; Verstra\"{e}te, J. Proof of a conjecture of Erd\H{o}s on triangles in set-systems. Combinatorica 25 (2005), no. 5, 599--614.


\bibitem{P2} Pikhurko, O. Exact Computation of the Hypergraph Tur\'{a}n Function for Expanded Complete 2-Graphs.
J. Combinatorial Theory  Ser. B 103 (2013) 220--225.


\end{thebibliography}
\end{document}